\documentclass[10pt]{article}

\usepackage[english]{babel}
\usepackage[latin1]{inputenc}
\usepackage{amsmath, latexsym}
\usepackage{amssymb}
\usepackage{rotating}

\usepackage[authoryear, round]{natbib}
\usepackage{color}
    \definecolor{myred}{rgb}{0.5,0,0}
    \definecolor{myblue}{rgb}{0,0,0.75}
    \definecolor{mygreen}{rgb}{0,0.5,0}
\newtheorem{theorem}{Theorem}[section]
\newtheorem{lemma}[theorem]{Lemma}
\newtheorem{remark}[theorem]{Remark}

\newtheorem{proposition}[theorem]{Proposition}

\newtheorem{assumption}[theorem]{Assumption}
\newtheorem{algorithm}[theorem]{Algorithm}

\numberwithin{equation}{section}

\newlength{\captionwidth}
\begin{document}

\title{Capital allocation for credit portfolios with kernel estimators}

\author{%
Dirk Tasche\thanks{Lloyds TSB Corporate Markets, Red Lion Court, 46-48 Park Street, London SE1 9EQ, United Kingdom.\newline
E-mail:
dirk.tasche@gmx.net
\newline The opinions expressed in this paper
are those of the author and do not necessarily reflect views of
Lloyds TSB.}}

\date{May 2008}
\maketitle

\begin{abstract}
Determining contributions by sub-portfolios or single exposures to
portfolio-wide economic capital for credit risk is an important risk measurement task.
Often economic capital is measured as Value-at-Risk (VaR) of the
portfolio loss distribution.
For many of the credit portfolio risk models used in practice,
the VaR contributions then have to be estimated from Monte Carlo
samples. In the context of a partly continuous loss distribution (i.e.\ continuous
except for
a positive point mass on zero), we investigate
how to combine kernel estimation methods with importance sampling
to achieve more efficient (i.e.~less volatile) estimation of VaR contributions.
\end{abstract}

%\renewcommand{\thefootnote}{\arabic{footnote}}
%\setcounter{footnote}{0}

%%%%%%%%%%%%%%%
% New section %
%%%%%%%%%%%%%%%
\section{Introduction}

In many financial institutions, there is a well established practice
of measuring the risk of their portfolios in terms of
\emph{economic capital} \citep[cf., e.g.][]{Dev2004}.
Measuring portfolio-wide economic capital, however, is only the
first step towards active, portfolio-oriented risk management.
For purposes like identification of risk concentrations, risk-sensitive
pricing or portfolio optimisation it is also necessary to decompose
portfolio-wide economic capital into a sum of risk contributions by
sub-portfolios or single exposures \citep[see, e.g.,][]{Litterman96}.

While already calculating or estimating economic capital is non-trivial in general,
determining risk decompositions is even more demanding \citep[see, e.g.,][]{Yamai02}.
For most of the economic capital models used in practice,
no closed-form solutions for risk contributions are available\footnote{%
See \citet{Tasche2004} or \citet{Tasche2006} for notable exceptions.
}.
As a consequence, for such models there is a need
for simulation to create samples from which to estimate economic capital as well as
risk contributions. These estimations often involve evaluations of very far tails
of the risk return distributions, causing high variability of the estimates. Various
variance reduction techniques have been proposed, one of the more popular being importance sampling
\citep[see][for its application to credit risk]{Glasserman&Li2005, Merino04,
Kalkbreneretal04}. \citet{Glasserman2005} suggests a two-step importance sampling approach to
the estimation of contributions to Value-at-Risk (VaR), the most popular metric
underlying economic capital methodologies. 

\citeauthor{Glasserman2005}'s approach,
however, does apply to discrete loss distributions only (i.e.~to distributions
such that each potential loss value has a positive probability to be assumed), which means a significant
restriction as loss distributions based on continuous loss given default rate distributions
seem to be more realistic in practical applications. 
According to Theorem \ref{pr:L} below such loss distributions have the property that each single
potential positive loss has probability zero to be assumed but there is a positive probability of not observing any loss.
The probability of not observing losses is significantly positive in particular for small and medium portfolio sizes (with 200 names 
or less) which occur, for instance, in typical securitisation deals.

To solve the problem of determining VaR contributions when positive losses have a
continuous distribution\footnote{%
Instead of trying to determine VaR contributions
in a continuous or semi-continuous setting, some risk managers use contributions to
expected shortfall. This approach is very fruitful and has some other
advantages \citep[see][]{Kalkbreneretal04, Merino04}.},
in this paper we follow, where possible, the path used by \citet{GL00}, \citet{EpperleinSmillie}, and \citet{GL06} who apply
kernel estimation methods for estimating VaR contributions and
contributions to spectral risk measures in a market risk context with
continuous distributions. {Kernel estimators are a well-established concept
\citep[see, e.g.,][]{Pagan99} to
deal with the issue of estimating non-elementary expectations as they occur
in the context of capital allocation.} Due to the
rare event issue characteristic for credit risk{, entailing rather volatile
estimates when standard Monte Carlo simulation is used}, we combine the kernel estimation
technique with importance sampling (shifting the means of the systematic factors to be more specific) for reducing estimation variance.
The paper is organized as follows:
\begin{itemize}
    \item Section \ref{se:capital} contains a review, in the necessary details, of the capital
    allocation problem  and the specific issues with the estimation
    of risk contributions when the loss distribution is partly continuous.
    \item Section \ref{se:kernel} provides a brief review of kernel estimators
    for densities and conditional expectations as well as its application
    to credit loss distributions with a partly continuous distribution.
    \item Section \ref{se:importance} introduces the model studied here
    and explains how to combine the kernel estimators from Section \ref{se:kernel}
    with importance sampling for credit risk.
    \item Application of the algorithms introduced is illustrated with
    a numerical example in Section \ref{se:example}.
    \item We conclude with  an assessment in Section \ref{se:concl} of what has been reached.
\end{itemize}

%%%%%%%%%%%%%%%
% New section %
%%%%%%%%%%%%%%%
\section{Capital allocation}
\label{se:capital}

In the following, we consider the following stochastic credit portfolio loss model:
\begin{equation}\label{eq:L}
    L\ =\ \sum_{i=1}^n L_i.
\end{equation}
$L_1, \ldots, L_n \ge 0$ are random variables that represent the losses
that a financial institution suffers on its exposures to borrowers
$i=1, \ldots, n$ within a fixed time-period, e.g.\ one year. The random variable
$L$ then expresses the portfolio-wide loss. We denote by $\mathrm{P}[\ldots]$
the real-world probability distribution that underlies model \eqref{eq:L}. In other
words, $\mathrm{P}[\ldots]$ is calibrated in such a way that it reflects
as close as possible observed loss frequencies. $\mathrm{P}[L \le \ell]$, for
instance, stands for the probability of observing portfolio-wide losses that do
not exceed the amount $\ell$. The operator $\mathrm{E_P}[\ldots]$ is defined as
mathematical expectation with respect to probability $\mathrm{P}$. In particular,
$\mathrm{E_P}[L]$ reflects the real-world probability weighted mean of the
portfolio-wide loss.

As mentioned in the introduction, it is common practice for financial institutions to measure
the risk inherent in their portfolios in terms of economic capital (EC). As credit risk,
for most institutions, is considered to be most important, this is in particular relevant for
credit portfolios. EC is commonly
understood as a capital buffer intended to cover the losses of the lending financial institution
with a high probability. This interpretation makes appear very natural the
definition
\begin{subequations}
\begin{equation}\label{eq:EC}
     \text{EC}\ = \ \mathrm{VaR}_{\mathrm{P},\alpha}(L) - \mathrm{E_P}[L],
\end{equation}
where the \emph{Value-at-Risk} (VaR) is given as a high-level (e.g.\ $\alpha = 99.9\%$)
quantile of the portfolio-wide loss:
\begin{equation}\label{eq:VaR}
    \mathrm{VaR}_{\mathrm{P},\alpha}(L)\ = \ \min\{\ell:\,\mathrm{P}[L\le\ell]\ge\alpha\}.
\end{equation}
\end{subequations}
Hence, if a financial institutions holds EC according to \eqref{eq:EC} and
charges the loans granted with upfront fees adding up to $\mathrm{E_P}[L]$, the probability
that it will lose all its EC is not higher than $1-\alpha$. Note that, despite its
intuitive appeal, VaR as a risk measure is criticised, e.g.\ for its potential lack
of rewarding diversification \citep[see][and the references therein]{Acerbi&Tasche}.

Active risk management involves more than just measuring portfolio-wide capital according
to \eqref{eq:EC}. Additionally, it is of interest to identify which
parts of the portfolio bind the largest portions of EC. The corresponding
process of determining a risk-sensitive decomposition of EC is called
\emph{capital allocation}. While for the expectation part $\mathrm{E_P}[L]$ of EC on
the right-hand side of \eqref{eq:EC} there is the natural decomposition
\begin{subequations}
\begin{equation}\label{eq:decom.E}
     \mathrm{E_P}[L]\ = \ \sum_{i=1}^n \mathrm{E_P}[L_i],
\end{equation}
there is no such obvious decomposition
\begin{equation}\label{eq:decom}
    \mathrm{VaR}_{\mathrm{P},\alpha}(L) \ = \ \sum_{i=1}^n \mathrm{VaR}_{\mathrm{P},\alpha}(L_i\,|\,L)
\end{equation}
\end{subequations}
for the VaR-part of EC into \emph{risk contributions}\footnote{%
\citet{Kalkbrener05} considers relation \eqref{eq:decom} in a more general context. He
calls it ``linear aggregation''.}. Indeed, the choice of the decomposition method depends
on the concept of risk sensitivity adopted. Interpreting risk sensitivity as
compatibility with portfolio optimization, \citet{Tasche1999} proved that the
risk contributions $\mathrm{VaR}_{\mathrm{P},\alpha}(L_i\,|\,L)$ on the right-hand side of
\eqref{eq:decom} should be defined as directional derivatives, i.e.
\begin{equation}\label{eq:contrib}
    \mathrm{VaR}_{\mathrm{P},\alpha}(L_i\,|\,L) \ = \
    \frac{d\, \mathrm{VaR}_{\mathrm{P},\alpha}(L + h\,L_i)}{d\, h}\bigl|_{h=0}.
\end{equation}
As VaR is a positively homogeneous\footnote{%
I.e.\ $\mathrm{VaR}_{\mathrm{P},\alpha}(h\,L) =
h\,\mathrm{VaR}_{\mathrm{P},\alpha}(L)$ for positive $h$.} risk measure, by Euler's theorem, then \eqref{eq:decom}
holds. \eqref{eq:contrib} displays the concept of risk contribution applied in this paper.
Note however that in general, and in particular if the distribution of $L$ has no density, the derivative on the right-hand side
\eqref{eq:contrib} need not exist. See, e.g., \citet[][Assumption (S)]{Tasche1999} for conditions ensuring existence of the derivative.
Depending on the objective of the portfolio analysis, other approaches to
determining risk contributions are reasonable, see, e.g., Section 3.1 of \citet{Tasche2006}
for an account of these.

In general, no closed-form representations of
$\mathrm{VaR}_{\mathrm{P},\alpha}(L)$ and  the risk contributions $\mathrm{VaR}_{\mathrm{P},\alpha}(L_i\,|\,L)$
are available. Therefore, often, these quantities can only be inferred from Monte-Carlo
samples. This means essentially to generate a sample
\begin{equation}\label{eq:MonteCarlo}
(L^{(t)}, L^{(t)}_1, \ldots, L^{(t)}_n),\quad t = 1, \ldots, T,
\end{equation}
and then to estimate the quantities under consideration on the basis of this sample.
How to do this is quite obvious for VaR, but is much less clear for the risk contributions
$\mathrm{VaR}_{\mathrm{P},\alpha}(L_i\,|\,L)$ as, in general, estimating derivatives of stochastic quantities without
closed-form representation is a subtle issue.

Fortunately, it turns out \citep{GL00, L99, Tasche1999}
that, under fairly general conditions on the joint distribution of
$L$ and $L_i$, the derivative \eqref{eq:contrib} coincides with an expectation of the
loss related to borrower $i$ conditional on the event of observing a portfolio-wide
loss equal to VaR.
\begin{equation}\label{eq:deriv}
    \frac{d\, \mathrm{VaR}_{\mathrm{P},\alpha}(L + h\,L_i)}{d\, h}\bigl|_{h=0} \ =\
    \mathrm{E_P}[L_i\,|\,L = \mathrm{VaR}_{\mathrm{P},\alpha}(L)]
\end{equation}
%Note that, when choosing the risk contribution of $L_i$ according to \eqref{eq:deriv}, the sum
%of the contributions will equal $\mathrm{VaR}_{\mathrm{P},\alpha}(L)$.
%
If $\mathrm{P}[L = \mathrm{VaR}_{\mathrm{P},\alpha}(L)]$ is positive, the conditional expectation on
the right-hand side of \eqref{eq:deriv} is given by
\begin{equation}\label{eq:cond_exp_elem}
    \mathrm{E_P}[L_i\,|\,L = \mathrm{VaR}_{\mathrm{P},\alpha}(L)]\ =\
    \frac{\mathrm{E_P}[L_i\,\mathbf{1}_{\{L = \mathrm{VaR}_{\mathrm{P},\alpha}(L)\}}]}{\mathrm{P}[L = \mathrm{VaR}_{\mathrm{P},\alpha}(L)]}.
\end{equation}
Even if $\mathrm{P}[L = \mathrm{VaR}_{\mathrm{P},\alpha}(L)]$ is positive, its magnitude will usually be very small,
such as $1-\alpha$ or less. \citet{Glasserman2005} shows how to apply importance sampling in such
a situation in order to efficiently estimate $\mathrm{E_P}[L_i\,|\,L = \mathrm{VaR}_{\mathrm{P},\alpha}(L)]$.

However, a crucial condition for \eqref{eq:deriv} to hold exactly is the existence of
a density of the distribution of $L$. The probability $\mathrm{P}[L = \mathrm{VaR}_{\mathrm{P},\alpha}(L)]$ then
equals zero, and consequently the right-hand side of \eqref{eq:cond_exp_elem} is undefined\footnote{%
This problem can be avoided by using the risk measure Expected Shortfall
\citep[see, e.g.,][]{Acerbi&Tasche} instead of VaR. With the definition of Expected Shortfall slightly simplified for practical purposes, \eqref{eq:cond_exp_elem}
then reads $\mathrm{E_P}[L_i\,|\,L \ge \mathrm{VaR}_{\mathrm{P},\alpha}(L)]\, =\,(1-\alpha)^{-1}
    \mathrm{E_P}[L_i\,\mathbf{1}_{\{L \ge \mathrm{VaR}_{\mathrm{P},\alpha}(L)\}}]$.}.
    In this situation, the conditional expectation $\mathrm{E_P}[L_i\,|\,L = \mathrm{VaR}_{\mathrm{P},\alpha}(L)]$
is still well-defined \citep[see, e.g., Remark 5.4 of][]{Tasche1999}, but its estimation from a sample like
\eqref{eq:MonteCarlo} requires more elaborated non-parametric methods. \citet{Mausser04}
suggest using an estimation method based on weighted combinations of order-statistics. We follow here
\citet{GL00} who applied kernel estimation
methods for VaR contributions when optimizing returns in a portfolio of stocks.
The kernel estimation procedures, however, have to be adapted to the rare-event
character of credit risk. Therefore,  in the remainder of the paper we modify the approach by
\citeauthor{GL00} in a way that
can be described as a combination of kernel estimation and importance sampling.

%%%%%%%%%%%%%%%
% New section %
%%%%%%%%%%%%%%%
\section{Kernel estimators}
\label{se:kernel}

In this section, we introduce the classical Rosenblatt-Parzen kernel estimator for densities
and the Nadaraya-Watson kernel estimator for conditional expectations in a way that links
naturally to the risk contribution concept of Section \ref{se:capital}. The general reference for
this section is \citet[][Chapters 2 and 3]{Pagan99}.

\subsection{The Rosenblatt-Parzen kernel estimator for densities}
\label{se:Rosenblatt}
Assume that $x_1, \ldots, x_T$ is
a sample of independent realisations of a random variable $X$ with density $f$. The Rosenblatt-Parzen estimator
$\hat{f}_h$ with bandwidth $h>0$ for $f$ can be constructed as follows:
\begin{itemize}
    \item Let $X^\ast$ be a random variable
    whose distribution is given by the empirical distribution corresponding to the sample $x_1, \ldots, x_T$,
    i.e.\ $\mathrm{P}[X^\ast = x_t] = 1/T$, $t = 1, \ldots, T$.
    \item Let $\xi$ a random variable with density (kernel) $\varphi$.
    \item Assume that $X^\ast$ and $\xi$ are independent.
    \item Then the estimator $\hat{f}_h$ is defined as the density of $X^\ast + h\,\xi$:
\begin{equation}\label{eq:dens}
        \hat{f}_h(x)\ = \hat{f}_{h, x_1, \ldots, x_T}(x)\ =\ \tfrac1{h\,T} \sum\nolimits_{t=1}^T \varphi\bigl(\tfrac{x-x_t}h\bigr).
\end{equation}
\end{itemize}
If $f$ and $\varphi$ are appropriately ``smooth'' \citep[see][Theorem 2.5 for details]{Pagan99},
it can be shown for
$h = h_T \xrightarrow{T\to\infty} 0$, $h_T\,T\xrightarrow{T\to\infty}\infty$ that $\hat{f}_{h_T}(x)$ is a pointwise mean-squared consistent
estimator of $f$, i.e.
\begin{equation}\label{eq:dens_consis}
     \lim_{T\to\infty}\mathrm{E}[(f(x) - \hat{f}_{h_T,X_1, \ldots, X_T}(x))^2]\, =\,  0,
     \quad x\in\mathbb{R},
\end{equation}
with independent copies $X_1, \ldots, X_T$ of $X$. While the Rosenblatt-Parzen density
estimator is rather robust with respect to the choice of
the kernel $\varphi$, it is quite sensitive to the choice of the bandwidth $h$. For the univariate case
we consider here, efficient techniques like cross validation for the choice of the bandwidth
are available. However, such more elaborated
techniques usually involve some optimisation procedures that can be very time-consuming for large
samples. As a consequence, for the purpose of this paper we confine ourselves to applying a simple
rule of thumb by Silverman \citep[cf.\ Chapter 2 in][]{Pagan99}
\begin{equation}\label{eq:Silverman}
    h\ =\ 1.06\,\sigma\,T^{-1/5},
\end{equation}%
where $\sigma$  denotes the standard deviation of the sample $x_1, \ldots, x_T$. Moreover, we choose the standard normal density
as the kernel $\varphi$.

\subsection{The Nadaraya-Watson kernel estimator for conditional expectations}
\label{se:Nadaraya}

Assume that $(x_1, y_1), \ldots, (x_T, y_T)$ is a sample of realisations of a random vector $(X,Y)$
where $X$ has a density $f$. The Nadaraya-Watson estimator $\hat{\mathrm{E}}_h[Y\,|\,X=x]$ with bandwidth
$h$ for $\mathrm{E}[Y\,|\,X=x]$ can be constructed as follows:
\begin{itemize}
    \item Let $(X^\ast, Y^\ast)$ a random vector whose distribution is given by the empirical distribution
    corresponding to the sample $(x_1, y_1), \ldots, (x_T, y_T)$, i.e.\
    $\mathrm{P}[(X^\ast, Y^\ast) = (x_t, y_t)] = 1/T$.
    \item Let $\xi$ a random variable with density (kernel) $\varphi$.
    \item Assume that $(X^\ast, Y^\ast)$ and $\xi$ are independent.
    \item Then the estimator $\hat{\mathrm{E}}_h[Y\,|\,X=x]$ is defined as the expectation
    of $Y^\ast$ conditional on $X^\ast + h\,\xi=x$:
    \begin{equation}\label{eq:est_cond_exp}
        \hat{\mathrm{E}}_h[Y\,|\,X=x]\ =\ \hat{\mathrm{E}}_{h, (x_1, y_1), \ldots, (x_T, y_T)}[Y\,|\,X=x]\ =\  \frac{\sum\nolimits_{t=1}^T y_t\,\varphi\bigl(\tfrac{x-x_t}h\bigr)}
                                          {\sum\nolimits_{t=1}^T \varphi\bigl(\tfrac{x-x_t}h\bigr)}.
\end{equation}
\end{itemize}
If $f$ and $\varphi$ are appropriately ``smooth'' \citep[see][Theorem 3.4 for details]{Pagan99},
it can be shown for
$h = h_T \xrightarrow{T\to\infty} 0$, $h_T\,T\xrightarrow{T\to\infty}\infty$, and $f(x) > 0$ that $\hat{\mathrm{E}}_h[Y\,|\,X=x]$ is a pointwise consistent
estimator of $\mathrm{E}[Y\,|\,X=x]$, i.e.
\begin{equation}\label{eq:exp_consis}
     \lim_{T\to\infty}\mathrm{P}\bigl[\,\bigr|\,\mathrm{E}[Y\,|\,X=x] - \hat{\mathrm{E}}_{h_T, (X_1, X_1), \ldots, (X_T, X_T)}[Y\,|\,X=x]\,\bigr| > \varepsilon\bigr]\, =\,  0,
     \quad \varepsilon >0\text{\ arbitrary},
\end{equation}
with independent copies $(X_1, Y_1), \ldots, (X_T, Y_T)$ of $(X,Y)$. The construction of the Nadaraya-Watson estimator
\eqref{eq:est_cond_exp} as described above allows to interpret the conditional expectation estimation problem
as an extended density estimation problem. This suggests to choose the same bandwidth $h$ and the same kernel
$\varphi$ for the estimators \eqref{eq:dens} and \eqref{eq:est_cond_exp}.

\begin{remark}\label{rm:sum}
Assume that in the random vector $(X, Y)$ the $X$-component is a sum of random
variables $X_1, \ldots, X_n$ and that we are interested in estimating
$\mathrm{E}[X_i\,|\,X=x]$, $i=1, \ldots, n$. Define $(X_1^\ast, \ldots, X_n^\ast)$,
analogously to $(X^\ast, Y^\ast)$ as the
``empirical'' version of $(X_1, \ldots, X_n)$. According to  \eqref{eq:est_cond_exp},
the Nadaraya-Watson estimator of $\mathrm{E}[X_i\,|\,X=x]$ can then be specified as
\begin{subequations}
\begin{equation}\label{eq:sum}
 \hat{\mathrm{E}}_h[X_i\,|\,X=x] \ = \ \mathrm{E}[X_i^\ast\,|\,X^\ast+h\,\xi=x],
\end{equation}
with an appropriate auxiliary variable $\xi$ independent of $(X_1^\ast, \ldots, X_n^\ast)$.
If the same bandwidth $h$ is applied for all $i$, then from representation \eqref{eq:sum} follows
\begin{equation}\label{eq:all_sum}
 \sum_{i=1}^n\hat{\mathrm{E}}_h[X_i\,|\,X=x] \ = \ x-h\,\mathrm{E}[\xi\,|\,X^\ast+h\,\xi=x].
\end{equation}
\end{subequations}
As we will note in Section \ref{se:comments} when commenting on Table \ref{tab:6}, the size of the difference of the left-hand side of \eqref{eq:all_sum} and $x$ can be regarded as providing additional information on
the choice of the bandwidth $h$. We will apply the multiplicative adjustment suggested by \citet[][Equation (7)]{EpperleinSmillie} to force additivity 
in the sense of \eqref{eq:decom} on the estimated contributions to VaR.
\end{remark}

\subsection{Application to credit losses}
\label{se:application}

For every credit risk portfolio, there is some positive probability of observing no losses. While for large portfolios,
this probability will usually be negligibly small, we will see by the example from Section \ref{se:example}  that
the probability of zero loss can be of significant magnitude for smaller portfolios. As a consequence, we cannot
assume that the loss variable $L$ from model \eqref{eq:L} has an unconditional density. For otherwise $\mathrm{P}[L = 0]$ would be
zero. Hence, at first glance, applying the estimators \eqref{eq:dens} and \eqref{eq:est_cond_exp} seems not possible
in the context of model \eqref{eq:L}. The following assumption, however, allows us to deal with this problem.

\begin{assumption} \label{as:L}
There is a random vector $(S_1, \ldots, S_k)$, (called \emph{systematic factors}), with the following
properties:
\begin{enumerate}
    \item[(i)] The loss variables $L_i$ in  \eqref{eq:L} are independent conditional on realisations of
    $(S_1, \ldots, S_k)$.
    \item[(ii)] For each $i=1, \ldots, n$, there are probabilities $p_i(s_1, \ldots, s_k) \in [0,1]$
    and densities $0 \le f_i(\ell, s_1, \ldots, s_k)$ such that
    the distribution function of $L_i$ conditional on
    $(S_1, \ldots, S_k)$ is given for $\ell \ge 0$ by
\begin{multline}\label{eq:cond_distr}
    \mathrm{P}[L_i \le \ell\,|\,(S_1, \ldots, S_k)=(s_1, \ldots, s_k)] = \\
    1- p_i(s_1, \ldots, s_k) + p_i(s_1, \ldots, s_k) \int_0^\ell f_i(x, s_1, \ldots, s_k)\,d x.
\end{multline}
\end{enumerate}
\end{assumption}
Under Assumption \ref{as:L}, it is easy to derive the following result on the representation of the
unconditional distribution of the portfolio-wide loss $L$.
\begin{theorem}\label{pr:L}
Define $I_n=\{1, \ldots, n\}$ and write $\bigotimes$ for the multiple convolution of densities. Then,
under assumption \ref{as:L}, for $\ell \ge 0$ the distribution function of the loss variable
$L$ from \eqref{eq:L} can be written as
\begin{align}\label{eq:rep_L}
    \mathrm{P}[L \le \ell]& =  p+(1-p) \int_0^\ell f(x)\,dx,\\
\intertext{with}
p & = \mathrm{P}[L_1=0, \ldots, L_n=0],\notag\\
f(x) & = \sum_{\emptyset\not= I\subset I_n} \mathrm{E_P}\Big[\prod_{i\in I}p_i(S_1, \ldots, S_k)
            \prod_{i\in I_n\backslash I}(1-p_i(S_1, \ldots, S_k))\, \bigl(\bigotimes_{i\in I} f_i(\cdot,S_1, \ldots, S_k)\bigr)(x)\Big]
            \notag.
\end{align}
\end{theorem}
Although, due to the involved multiple convolutions, \eqref{eq:rep_L} is not really useful
for calculating the distribution of $L$, it allows us to assume that, conditional on being positive,
the portfolio-wide loss has a density, i.e.\ for $\ell \ge 0$
\begin{equation}\label{eq:cond_dens_L}
    \mathrm{P}[L \le \ell\,|\,L >0]\ =\ \int_0^\ell f(x)\,dx.
\end{equation}
Define $\mathrm{P}^\ast$ by
\begin{equation}\label{eq:Pstar}
\mathrm{P}^\ast[A]\ =\ \mathrm{P}[A\,|\,L >0]
\end{equation}
for any relevant event $A$. The following lemma is then obvious.

\begin{lemma}\label{le:Pstar}
If the probability $\mathrm{P}^\ast$ is given by \eqref{eq:Pstar}, then for $\ell > 0$ the expectations
conditional on $L = \ell$ with respect to $\mathrm{P}^\ast$ and $\mathrm{P}$ are identical. In particular,
in the context of model \eqref{eq:L} for all $i = 1, \ldots, n$ and $\ell > 0$ we have
\begin{equation*}
        \mathrm{E_P}[L_i\,|\,L = \ell]\ = \ \mathrm{E}_{\mathrm{P}^\ast}[L_i\,|\,L = \ell].
\end{equation*}
\end{lemma}
Note that a sample \eqref{eq:MonteCarlo} generated under measure $\mathrm{P}$ becomes a sample generated
under measure $\mathrm{P}^\ast$ when all $(n+1)$-tuples $(L^{(t)}, L_1^{(t)}, \ldots, L_n^{(t)})$ with
$L^{(t)} = 0$ are eliminated. For this sub-sample then the preconditions for applying estimators
\eqref{eq:dens} and \eqref{eq:est_cond_exp} are satisfied. This observation leads to the following algorithm
for estimating the risk contributions according to \eqref{eq:deriv} for model \eqref{eq:L} by Monte Carlo
sampling.
\begin{algorithm}\label{al:standard}\
\vspace{-2ex}
\begin{enumerate}
    \item Generate a sample like \eqref{eq:MonteCarlo} from the real-world probability measure $\mathrm{P}$.
    \item Determine an estimate\footnote{%
    We may assume $\hat{\ell} > 0$ as for real-world portfolios the case $\hat{\ell} = 0$
    seems very unlikely.} $\hat{\ell}$ of $\mathrm{VaR}_{\mathrm{P}, \alpha}(L)$ from this sample.
    \item Extract the sub-sample with $L^{(t)} > 0$ from the previous sample.
    \item Calculate, on the basis of the sub-sample with $L^{(t)} > 0$, the bandwidth $h^\ast$ according
    to \eqref{eq:Silverman}.
    \item Calculate the estimates for $\mathrm{E_P}[L_i\,|\,L = \mathrm{VaR}_{\mathrm{P},\alpha}(L)]$, $i = 1, \ldots, n$,
    on the basis of the sub-sample with $L^{(t)} > 0$, according to \eqref{eq:est_cond_exp} as
    $\hat{\mathrm{E}}_{h^\ast}[L_i\,|\,L=\hat{\ell}]$ (cf.~Remark \ref{rm:sum}).
\end{enumerate}
\end{algorithm}
In the following section, we will modify this algorithm by incorporating importance sampling for
reducing the variances of the estimates.

%%%%%%%%%%%%%%%
% New section %
%%%%%%%%%%%%%%%
\section{Importance sampling for credit risk}
\label{se:importance}

\citet{McNeil05}, at the beginning of Chapter 8.5, write ``A
possible method for calculating risk measures and related
quantities such as capital allocations is to use Monte Carlo (MC)
simulation, although the problem of \emph{rare event simulation}
arises.'' They then explain, in the context of risk contributions
to Expected Shortfall, that ``the standard MC estimator $\ldots$
will be unstable and subject to high variability, unless the
number of simulations is very large. The problem is of course that
most simulations are `wasted', in that they lead to a value of $L$
which is smaller than $\mathrm{VaR}_{\mathrm{P}, \alpha}(L)$.''
What applies to contributions to expected shortfall applies even
more to VaR contributions as these, according to \eqref{eq:deriv},
are related to events still more rare and actually in our case of probability
zero. Thus, the key idea with importance sampling is to replace
the real-world probability measure $\mathrm{P}$ by a probability measure
$\mathrm{Q}$ which puts more mass on the interesting events.

Some technical
assumptions are needed to guarantee that such a replacement of
probabilities does really work.

\begin{assumption}\label{as:imp_sampling} \
\vspace{-1ex}
\begin{itemize}\setlength{\itemsep}{0.5ex}
    \item The probability measures $\mathrm{P}$ and $\mathrm{Q}$
are defined on the same measurable space $(\Omega, \mathcal{F})$.
    \item There is a measure $\mu$ on $(\Omega, \mathcal{F})$ such that both
$\mathrm{P}$ and $\mathrm{Q}$ are absolutely continuous with
respect to $\mu$. Denote by $f$ the density of $\mathrm{P}$ and by
$g$ the density of $\mathrm{Q}$.
    \item $f >0$ implies $g> 0$, i.e.\ $\mathrm{P}$ is absolutely
    continuous with respect to $\mathrm{Q}$.
\end{itemize}
\end{assumption}
Under Assumption \ref{as:imp_sampling}, the \emph{likelihood
ratio}
\begin{equation}\label{eq:R}
    R \ =\ \frac f g
\end{equation}
is $\mathrm{Q}$-almost surely well-defined. This
implies that any expectation with respect to $\mathrm{P}$ can be
expressed as an expectation with respect to $\mathrm{Q}$, i.e.\
for any integrable $X$ holds
\begin{equation}\label{eq:replace}
    \mathrm{E_P}[X] \ = \ \mathrm{E_Q}[R\,X].
\end{equation}
However, according to \eqref{eq:deriv}, for the purpose of this paper
it is a \emph{conditional}
rather than an unconditional expectation we are interested in. The following
proposition gives a result analogous to \eqref{eq:replace}, for conditional
expectations.

\begin{proposition}\label{pr:cond_exp}
Let $\mathrm{P}$ and $\mathrm{Q}$ be probability measures as in Assumption \ref{as:imp_sampling},
with $\mu$-densities $f$ and $g$ respectively. Define the likelihood ratio $R$ by \eqref{eq:R}.
Then we have for any sub-$\sigma$-algebra
$\mathcal{A}$ of $\mathcal{F}$ and any integrable random variable $X$
\begin{equation}\label{eq:cond_exp}
    \mathrm{E_P}[X\,|\,\mathcal{A}]\ = \ \frac{\mathrm{E_Q}[R\,X\,|\,\mathcal{A}]}{\mathrm{E_Q}[R\,|\,\mathcal{A}]}.
\end{equation}
\end{proposition}
\textbf{Proof.} See e.g.\ \citet[][Theorem 10.8]{Klebaner}.\hfill $\Box$

The choice of an appropriate measure $\mathrm{Q}$ for use in \eqref{eq:replace} and \eqref{eq:cond_exp}
is not at all obvious. 
Inspecting \eqref{eq:L} and \eqref{eq:deriv}, it becomes clear that, for the purpose of this paper, we are interested
in events that result in losses close\footnote{%
In general, $\hat{\ell}$ itself will have to be estimated. Thus, at
first glance, it seems strange to choose it as a basis for finding
$\mathrm{Q}$. However, in a first step, for instance, it can be replaced by a
rough estimate and be refined in further stages of the estimation
procedure.} to $\hat{\ell}=\mathrm{VaR}_{\mathrm{P},\alpha}(L)$, as defined by \eqref{eq:VaR}.

In Section \ref{se:example}, we will specifically consider a numerical example that
satisfies the following conditions which reflect assumptions commonly made for
models {in} industry.

\begin{assumption}\label{as:dep_tilting}
The portfolio-wide loss $L$ is given by \eqref{eq:L}. $L_1, \ldots, L_n$ are given as
\begin{equation*}
    L_i\ = \ A_i\,\mathbf{1}_{D_i}.
\end{equation*}
$D_1, \ldots, D_n$ are independent (default) events, \emph{conditional on a set of systematic factors}
$(S_1, \ldots, S_k)$, with $\mathrm{P}[D_i] = p_i$. The loss severity
variables $A_1, \ldots, A_n$ are positive and independent, as well as independent of the $D_1, \ldots,
D_n$ and the $S_1, \ldots, S_k$. The distribution of $A_i$ is specified via its density $a_i(s) \ge 0$.
\end{assumption}
Note that,
under the condition of the $A_i$ having densities, Assumption \ref{as:dep_tilting} is
a special case of Assumption \ref{as:L}. As a consequence of this, according to Theorem \ref{pr:L} the loss variable
$L$ from \eqref{eq:L} has, in theory, a density for its positive realisations.

\citet{Merino04} and \citet{Glasserman&Li2005} suggest a nested simulation procedure where exponential
twisting is applied for the estimation of expectations conditional on the systematic factors.%
The resulting conditional loss distribution then can essentially be described again by
Assumption \ref{as:dep_tilting}, with independence instead of conditional independence and modified probabilities of
default. Unfortunately, this nested
procedure cannot be applied for our problem of estimating $\mathrm{E_P}[L_i\,|\,L=\cdot]$ as there is
no independence conditional on $L$. In general, we have $\sigma(L) \not\subset \sigma(S_1, \ldots, S_k)$
and $ \sigma(S_1, \ldots, S_k) \not\subset \sigma(L)$. Therefore, no nesting of conditioning is
applicable either. Exponential twisting can be applied nevertheless but does not yield satisfactory results
as it is not clear how an optimal tilting parameter should be determined.

Approaches by
\citet{Kalkbreneretal04} and \citet{Glasserman&Li2005} are promising alternatives to exponential
twisting.
In these approaches,  the importance sampling measure
$\mathrm{Q}$ is created by changing the means of the systematic factors from Assumption \ref{as:dep_tilting}.
\citet{Kalkbreneretal04} suggest to determine the factor means appropriate for importance sampling by
solving a minimisation problem for the variance of the estimator they consider.  \citet{Glasserman&Li2005} instead
look for factor means that make the mode of the factor distribution coincide (approximately) with the mode of the
``zero-variance importance sampling'' distribution.
However, these approaches are  more complex compared to exponential
twisting as they involve -- in the case of a multi-factor model --
choosing several parameters  instead of only one as required by exponential
twisting.

The approach we follow in this paper is close in spirit to the
approach by \citet{Glasserman&Li2005}. 
\begin{subequations}
Define the importance sampling probability measure $\mathrm{Q} = \mathrm{Q}_\mu$ by Assumption \ref{as:dep_tilting},
but replace the vector $(S_1, \ldots, S_k)$ of systematic factors by a vector $(S_1^\ast, \ldots, S_k^\ast)$ with
\begin{equation}\label{eq:vector}
    S_i^\ast \ =\ S_i - \mathrm{E}_{\mathrm{P}}[S_i] + \mu_i
\end{equation}
where $\mu = (\mu_1, \ldots, \mu_k)$ satisfies
\begin{equation}\label{eq:mu}
    \mu_i \ \approx\ \mathrm{E}_{\mathrm{P}}[S_i\,|\,L = \mathrm{VaR}_{\mathrm{P}, \alpha}(L)].
\end{equation}
\end{subequations}
The conditional expectations in \eqref{eq:mu} are estimated by applying the Nadaraya-Watson estimator
\eqref{eq:est_cond_exp}.

Algorithm \ref{al:standard} for estimating the risk contributions according to \eqref{eq:deriv}
for model \eqref{eq:L} by Monte Carlo
sampling has to be modified as follows when importance sampling is applied.
\begin{algorithm}\label{al:importance}\
\vspace{-2ex}
\begin{enumerate}
    \item Generate a sample $(L^{(t)}, L^{(t)}_1, \ldots, L^{(t)}_n)$, $t = 1, \ldots, T_1$,
 from the original sampling probability measure $\mathrm{P}$.
    \item Determine an estimate
    $\hat{\ell}$ of $\mathrm{VaR}_{\mathrm{P}, \alpha}(L)$ from this sample.
    \item Calculate, on the basis of the sub-sample with $L^{(t)} > 0$, the standard MC simulation bandwidth $h^\ast$ according
    to \eqref{eq:Silverman}.
    \item Estimate the mean shift parameters $\mu_i$ according to \eqref{eq:mu} from the sample
    $(L^{(t)}, L^{(t)}_1, \ldots, L^{(t)}_n)$, $t = 1, \ldots, T_1$, for instance by applying the Nadaraya-Watson
    estimator.
    \item Generate a sample\footnote{%
    The new sample need not necessarily be generated by a new Monte Carlo simulation. Alternatively,
    as done in Section \ref{se:example}, the new sample can be created from the previous sample by substituting $(S_1^\ast, \ldots, S_k^\ast)$ from
    \eqref{eq:vector} for $(S_1, \ldots, S_k)$. In this case $T_2 = T_1$.%
    } $(L^{(t)}, L^{(t)}_1, \ldots, L^{(t)}_n, R^{(t)})$, $t = 1, \ldots, T_2$,
 from the importance sampling probability measure $\mathrm{Q}_\mu$. $R^{(t)}$ denotes realisations
 of the likelihood ratio as defined by \eqref{eq:R}.
    \item (Optional) Determine a refined estimate\footnote{%
    For instance, by ordering the pairs $(L^{(t)}, R^{(t)})$ in descending
    order according to the $L$-component and selecting the largest $t$ with
    $\sum_{k=1}^t R^{(k)}\le 1-\alpha$. Take then $L^{(t)}$ as the estimator.}
    $\hat{\ell}$ of $\mathrm{VaR}_{\mathrm{P}, \alpha}(L)$ from this sample.
    %
%    \item Extract the sub-sample with $L^{(t)} > 0$ from the sample.
    %
    \item Calculate, on the basis of the sub-sample with $L^{(t)} > 0$, the importance sampling bandwidth $h^\ast$ according
    to \eqref{eq:Silverman}.
    \item Estimate, on the basis of the sub-sample with $L^{(t)} > 0$, the $\mathrm{Q}_\mu$-conditional
    expectations
    $\mathrm{E}_{\mathrm{Q}_\mu}[R\,|\,L = \hat{\ell}\,]$ and
    $\mathrm{E}_{\mathrm{Q}_\mu}[R\,L_i\,|\,L = \hat{\ell}\,]$, $i = 1, \ldots, n$ according to
        \eqref{eq:est_cond_exp} with bandwidth $h^\ast$.
    \item Calculate the estimates for $\mathrm{E_P}[L_i\,|\,L = \mathrm{VaR}_{\mathrm{P},\alpha}(L)]$, $i = 1, \ldots, n$
     according to \eqref{eq:cond_exp}, inserting the estimates from the previous step.
\end{enumerate}
\end{algorithm}

%%%%%%%%%%%%%%%
% New section %
%%%%%%%%%%%%%%%
\section{Numerical example}
\label{se:example}

This section is divided into two parts. In Sub-section \ref{se:description} we provide
details of the portfolio model that is used for the simulation study and describe the
numerical results we seek to obtain. In Sub-section \ref{se:comments} we comment on the results
and point out their essential features.

\subsection{Description of simulation study}
\label{se:description}

We examine the performance of importance sampling estimation of VaR contributions
under Assumption \ref{as:dep_tilting}. Due to restriction in computational power,
the portfolio we consider is relatively small, with 96 assets. Note that, with such a portfolio
size, there may be a significant positive probability of not observing any losses, demonstrating that 
the issue tackled in Lemma \ref{le:Pstar} has some relevance.

The dependence structure of the portfolio is determined by four correlated systematic factors
$(S_1, S_2, S_3, S_4)$ which are each standard normal and
are jointly normally distributed with correlation matrix
\begin{equation}\label{eq:corr_matrix}
    \begin{pmatrix}
        1 & 0.75 & 0.05 & 0.05 \\
        0.75 & 1 & 0.05 & 0.05 \\
        0.05 & 0.05 & 1 & 0.25\\
        0.05 & 0.05 & 0.25 & 1
    \end{pmatrix}.
\end{equation}
According to \eqref{eq:corr_matrix}, factors $S_1$ and $S_2$ are strongly correlated,
factors $S_3$ and $S_4$ are moderately correlated, and the pairs $(S_1, S_2)$ and $(S_3, S_4)$
are weakly dependent.

Each factor corresponds to one sector that includes 24 assets. In each of the four sectors,
the risk characteristics of the 24 assets are specified as shown in Table \ref{tab:1}.
According to Table \ref{tab:1}, each of the four sectors in the portfolio
includes 12 high (2\%) PD \emph{(probability of default)} assets and 12 low (0.5\%) PD assets. In both of these two sub-sectors
there are 4 high exposure (25\$), 4 medium-size exposure (5\$), and 4 low exposure (1\$) assets.
Of the 4 high exposure assets, all have equal LGD \emph{(loss given default)} mean 50\% but 2 assets have high LGD variance
(12.5\%) and 2 assets have low
LGD variance (3.125\%). Similarly, among the 4 medium-size (low) exposure assets of equal LGD
mean, there are
2 high LGD variance and 2 low LGD variance assets. Note that for each combination of sector, PD,
exposure, and LGD variance there are two assets with identical risk characteristics.
This feature of the portfolio composition is intended to deliver a rough assessment of estimation
uncertainty due to the Monte Carlo simulation. For, by symmetry,
risk contributions for assets with identical risk characteristics should be equal but will not be
when estimated by Monte Carlo simulation.

\begin{assumption}\label{as:port}
For any asset $i \in \{1, 2, \ldots, 96\}$ in the portfolio, its individual loss variable $L_i$
in the sense of Assumption \ref{as:dep_tilting} is specified by the asset's sector $k(i) \in \{1,2,3,4\}$,
probability of default $PD_i$, exposure $v_i$, mean loss given default $LGD_i$,
and LGD variance $varLGD_i$ (as shown in Table
\ref{tab:1}) as follows:
\begin{itemize}
  \item The default event of the asset is given by $D_i = \{ \sqrt{r}\,S_{k(i)} + \sqrt{1-r}\,\xi_i \le
    \Phi^{-1}(PD_i)\}$ where $\Phi$ denotes the standard normal distribution function, $r = 0.18$ (equal for
    all $i$), and $\xi_1, \xi_2, \ldots, \xi_{96}$ are i.i.d.\ standard normal.
  \item The loss severity variable is given by  $A_i = v_i\,B_i$ where $B_1, B_2, \ldots, B_{96}$ are
  independent beta-distributed random variables with $\mathrm{E}[B_i] = LGD_i$ and $\mathrm{var}[B_i] =
  varLGD_i$.
\end{itemize}
\end{assumption}
The constant $r$ in the definition of the default events is the loading of the systematic risk in this model.
Its value was chosen with a view on
the asset correlations in the Basel II corporate risk weight formula which have a range from 0.12 to 0.24 \citep[][paragraph 272]{BaselAccord}.
The loss
severity distributions of the high LGD variance assets in Assumption \ref{as:port} are U-shaped (beta parameters
$a = 0.5$ and $b= 0.5$), the severity distributions of the low LGD variance assets are bell-shaped (beta parameters
$a = 3.5$ and $b= 3.5$).

By portfolio construction, on the one hand Sectors 1 and 2 are identical with respect
to their composition and risk characteristics. On the other hand, this holds also for Sectors 3 and 4.
As Sectors 1 and 2 are, however, stronger correlated than Sectors 3 and 4, one should expect that due to this  concentration the risk
contributions of assets in Sectors 1 or 2 are higher than the risk contributions of the corresponding assets
in Sectors 3 or 4.

The simulation exercise\footnote{%
The R-script for the calculations can be down-loaded at \texttt{http://www-m4.ma.tum.de/pers/tasche/}.%
} we conduct is structured as follows:
\begin{enumerate}
	\item We run 25 times a Monte Carlo simulation with 50,000 joint realisations of the losses of all
	the assets in the portfolio.
	\item In each simulation run, the following quantities are estimated:
\begin{itemize}
	\item Portfolio-wide VaR at 99.9\% level (by standard MC and importance sampling MC, see Table \ref{tab:2} for the results).
	\item For each factor, the conditional expectation $\mathrm{E}_{\mathrm{P}}[S_i\,|\,L = \mathrm{VaR}_{\mathrm{P}, \alpha}(L)]$ of the factor conditional on $L = \mathrm{VaR}_{\mathrm{P}, \alpha}(L)$ (by standard MC, Table \ref{tab:3}).
    \item For each sector, sector-stand-alone VaR at 99.9\% level (by standard MC and importance sampling MC, Table \ref{tab:5}).
    \item For each asset, the contribution of the asset to portfolio-wide VaR at 99.9\% level (by standard MC and importance sampling MC, Table \ref{tab:7}). Additionally, for each asset, the contribution of the asset to portfolio-wide VaR at 99.9\% level by importance sampling MC with reduced and enlarged respectively bandwidth for the conditional expectation, Table \ref{tab:8}).
\end{itemize}
\item Based on the estimates of step 2, in each simulation run the following quantities are
calculated:
\begin{itemize}
    \item The portfolio-wide mean of the loss, standard deviation of the loss, and the probability of observing a loss\footnote{%
For further reduction of the estimation variance, the standard deviation and sample size of the positive losses for the bandwidth (both for standard MC as well as for importance sampling) according to \eqref{eq:Silverman} are not estimated but numerically calculated. Under Assumption \ref{as:port} this can be done exactly for the standard deviation of the positive losses and approximately for the sample size of the positive losses, by approximating the distribution of the number of defaults via moment matching by a negative binomial distribution.%
} (Table \ref{tab:4}).
    \item The ratio of the sum of the VaR contributions of all assets and portfolio-wide VaR (Table \ref{tab:6}).
    \item For each sector, the contribution of the sector to portfolio-wide VaR at 99.9\% level (as the sum of the VaR contributions
  of the assets in the sector, Table \ref{tab:9}).
  \item The portfolio-wide diversification index\footnote{%
  See \citet[][Section 4]{Tasche2006} for a motivation of this definition and the definition of marginal diversification indices. The diversification indices are calculated on an \emph{unexpected loss (UL)} basis, in accordance with definition \eqref{eq:EC} of economic capital. Note that the value of the portfolio-wide diversification index depends upon whether the portfolio is decomposed into assets or into sectors since in the latter case the diversification potential is larger.} (Table \ref{tab:10})
\begin{subequations}
  \begin{equation}\label{eq:div_index}
    \mathrm{I}(L)\ =\ \frac{\mathrm{VaR}_{\mathrm{P}, \alpha}(L)-\mathrm{E}[L]}{\sum_{i=1}^n (\mathrm{VaR}_{\mathrm{P}, \alpha}(L_i)-\mathrm{E}[L_i])}.
  \end{equation}
  \item For each asset $i$, the marginal diversification index (Table \ref{tab:11})
  \begin{equation}\label{eq:div_index_asset}
    \mathrm{I}(L_i\,|\,L)\ =\ \frac{\mathrm{VaR}_{\mathrm{P},\alpha}(L_i\,|\,L)-\mathrm{E}[L_i]}{\mathrm{VaR}_{\mathrm{P}, \alpha}(L_i)-\mathrm{E}[L_i]}.
  \end{equation}
\end{subequations}
  \item For each sector, its marginal diversification index (defined as the sum of the VaR contributions
  of the assets in the sector divided by the VaR of the sector, Table \ref{tab:12}).
\end{itemize}
Marginal diversification indices as defined by \eqref{eq:div_index_asset} represent a direct application of 
risk contributions for risk concentration analysis. By construction, $\mathrm{I}(L_i\,|\,L) > \mathrm{I}(L)$ implies 
that reduction of the exposure to asset $i$ will improve portfolio diversification \citep[][Section 4]{Tasche2006}. Here, ``diversification''
is understood in a relative sense, namely comparing the actual economic capital assigned to the portfolio under consideration 
to the economic capital assigned to a worst case portfolio composed of co-monotonic loss variables.
\item We report three results for each of the values that is estimated or calculated from estimates:
\begin{itemize}
  \item The result of the first simulation run.
  \item The mean of the results of all 25 simulation runs (as approximation of the true value).
  \item The coefficients of variation\footnote{%
  The coefficient of variation of a sample is defined as the ratio of the sample standard deviation
  and the sample mean.%
   } of the results of all 25 simulation runs (as measure of estimation
  uncertainty).
\end{itemize}
\end{enumerate}
\subsection{Comments on the results}
\label{se:comments}

\paragraph{Table \ref{tab:2}.} With respect to the estimation of portfolio-wide VaR, according to the results displayed in Table \ref{tab:2}, there is a minor advantage in using importance sampling. With multi-step importance sampling as suggested for instance by \citet{Glasserman&Li2005} this advantage could be increased. However, the purpose of this paper is to deal with the estimation of VaR contributions. Therefore, here we need not look into the details of efficient estimation of VaR itself. As indicated by the relatively small sample variation coefficients both for standard as well as importance Monte Carlo sampling, the estimates in the first simulation run do not differ much from the means of 25 simulation runs.
\paragraph{Table \ref{tab:3}.} Table \ref{tab:3} shows that the estimates for the shift constants of the distribution of the systematic factors according to \eqref{eq:vector} and \eqref{eq:mu} are not very stable. This is indicated both by the high sample variation coefficients as well as by the 1st run estimate of the constant for the fourth factor which differs much from the sample mean. As the strongly correlated factors 1 and 2 might contribute more to portfolio risk as measured by VaR than the factors 3 and 4, it is no surprise that the indicated shifts of factors 1 and 2 are larger
than the ones of factors 3 and 4.
\paragraph{Table \ref{tab:4}.} The difference between the loss distribution under the original probability measure $\mathrm{P}$ and the importance sampling measures is demonstrated by Table \ref{tab:4}. The mean loss grows by more than ten times and then slightly overshoots the portfolio VaR under the original probability measure $\mathrm{P}$. The loss standard deviation increases by four times. The proportion of portfolio loss realisations with positive losses grows from 60\% to almost 100\%. Hence the proportion of the sample that can be used for kernel estimation is much greater in the case of importance sampling. Note that the characteristics of the shift loss distribution do not seem to vary much in the 25 simulation runs.
\paragraph{Table \ref{tab:5}.} By construction, from a risk perspective the four portfolio sectors are identical when considered stand-alone. This is confirmed by the estimates of stand-alone sector VaR  as displayed in Table \ref{tab:5}. The significantly higher sample variation coefficients of sectors 3 and 4 in the case of importance sampling show that the shift factor distribution according to \eqref{eq:vector} and \eqref{eq:mu} is not very well suited for the estimation of stand-alone sector characteristics.
\paragraph{Table \ref{tab:6}.} As indicated in Remark \ref{rm:sum}, although in theory the sum of the VaR contributions according to \eqref{eq:deriv} equals VaR, a sum of VaR contribution estimates made by kernel estimation can differ from VaR. Table \ref{tab:6} displays for different estimation approaches how large the difference can be. In general, it is larger for importance sampling and increases with the kernel estimation bandwidth. As we will see, the choice of the importance sampling bandwidth according to Silverman's rule of thumb \eqref{eq:Silverman} seems to be a reasonable compromise between reduction of sample variation by oversmoothing and unbiasedness as measured by the difference between the sum of VaR contributions and VaR. Note that, to make comparable the VaR contribution estimates by different approaches, the contributions as displayed in Tables \ref{tab:7}, \ref{tab:8}, and \ref{tab:9} are normalised such that their sum equals portfolio VaR.
\paragraph{Tables \ref{tab:7} and \ref{tab:8}.} As by symmetry results for sectors 2 and 4 are not essentially different from the results for sectors 1 and 3 respectively, Tables \ref{tab:7} and \ref{tab:8} display the results of VaR contribution estimates at asset level only for sectors 1 and 3. The tables allow to compare the estimation performance as yielded by four different approaches: kernel estimation based on standard Monte Carlo sampling, kernel estimation based on importance sampling with bandwidth chosen according to Silverman's rule of thumb, kernel estimation based on importance sampling with reduced bandwidth, and kernel estimation based on importance sampling with enlarged bandwidth. The average results of all the approaches do not differ too much -- which is also a consequence of the normalising applied to have the sum of the contributions equal to portfolio VaR. The variation of the estimates, as indicated by their coefficients of variation, is clearly highest for standard Monte Carlo, followed by importance sampling with 50\% of Silverman's bandwidth. Importance sampling with 200\% of Silverman's bandwidth has lower sample variation than importance sampling with 100\% of Silverman's bandwidth. As noticed above, the latter, however, need not be so much adjusted as the former. For this reason, importance sampling with bandwidth according to Silverman's rule might be considered the best estimation method for VaR contributions among the four approaches discussed here.

Taking into account that the values in Tables \ref{tab:7} and \ref{tab:8} were calculated on the basis of 25 simulation runs each with 50,000 loss realisations, the results are somewhat disappointing as the sample variation is still large in particular for small exposures with low probability of default. As demonstrated by Figure~\ref{fig:1}, there is nevertheless a clear gain in estimation efficiency by the application of importance sampling. Additionally, in contrast to the standard Monte Carlo approach, the importance sampling approaches have no problem with seeming zero VaR contributions of some assets (compare the first run results).

Note that the correlation structure of the portfolio is clearly reflected in the VaR contributions as the contributions by assets in sector 1 are significantly higher than the contributions by the assets with similar risk characteristics in sector 3. Note also that LGD variance has a strong impact on an asset's VaR contribution, in particular for those assets with the high exposures. The greater the LGD variance, the greater the VaR contribution.
\paragraph{Table \ref{tab:9}.} Table \ref{tab:9} displays at sector level what Table \ref{tab:7} shows at asset level. Compared to the asset level, sample variation at sector level is slightly less. Again, by importance sampling there is some gain in estimation efficiency compared to standard sampling. Note in particular for sector 1 the misleading first run standard estimate of the VaR contribution which seems to indicate that sector 1 is less correlated to the rest of the portfolio than sectors 2, 3, or 4.
\paragraph{Table \ref{tab:10}.} As demonstrated in Table \ref{tab:10}, estimates of portfolio-wide diversification indices in asset and sector context are fairly stable, and there is not much difference in efficiency between standard and importance sampling. The specification of the context in which the diversification indices are calculated indicates the scope of possible actions for a reduction of portfolio concentration. Sector context means that only the relative weights of the sectors may be changed but not relative weights of single assets. Superficially, the sector diversification index looks worse. This is caused by the fact that its denominator -- compared to the asset context -- is less because the sector VaR figures already incorporate a lot of diversification. In general, it does not make sense to compare diversification indices that were calculated in different contexts. Portfolio-wide diversification indices should rather be compared to the corresponding asset or sector diversification indices because this way guidance can be provided on how to change the portfolio for better diversification.
\paragraph{Table \ref{tab:11}.} Not surprisingly, as Table \ref{tab:11} shows importance sampling estimates are more efficient than standard sampling estimates also for the estimation of asset-level marginal diversification indices. In particular, importance sampling avoids observing negative diversification indices even in the first run estimates where estimation is much more uncertain than in the ``mean of all runs'' columns. Note however, that also the importance sampling first run estimates are misleading in so far as they seem to indicate that asset 4 (with low LGD variance) in sector 1 is the most dangerous in the portfolio. In fact, assets 1 and 2 in this sector with same exposure size, PD, and mean LGD are more dangerous -- as is correctly shown in columns 4 and 7 of Table \ref{tab:11} -- because they have got higher LGD variance. Note also that exposure concentration can ``override'' concentration caused by correlation. This is illustrated by assets 1 to 4 of sector 3 whose diversification indices are higher than the portfolio-wide index although sector 3 clearly has less correlation to the rest of the portfolio than sectors 1 or 2.
\paragraph{Table \ref{tab:12}.} Also at sector level the importance sampling estimates of the marginal diversification indices display less variation than the standard sampling estimates. In particular, the standard sampling first run estimate of the index for sector 1 is lower than the portfolio-wide diversification index. The misleading conclusion could be that shifting weight to sector 1 contributes to portfolio diversification, in spite of the strong correlation between sectors 1 and 2. In contrast, the importance sampling first run estimate of the diversification index for sector 1 correctly indicates that the index is larger than the portfolio-wide index. However, it also indicates erroneously that sector 1 is worse for portfolio diversification than sector 2.

%%%%%%%%%%%%%%%
% New section %
%%%%%%%%%%%%%%%
\section{Conclusions}
\label{se:concl}

In general, determining VaR contributions in a credit portfolio risk model that involves continuous
loss given default rate distributions is a non-trivial task. In the context of the
common approach by means of Monte-Carlo simulation, we have discussed how to adapt kernel
estimation methods for this problem and how to combine them with importance sampling.
Importance sampling, in the form of a shift
of the distribution of the systematic factors,
is applied here since the variability of the estimates is quite strong,
as a consequence of the rare-event character of credit risk realisations.

The numerical example presented in Section \ref{se:example} illustrates that the gain
in estimation efficiency by these methods is significant. It also reveals, however, that
the results yielded with these methods are not yet too satisfactory in so far
as in particular the variability of estimates of VaR contributions for exposures
with very small PDs remains still quite large. It seems worthwhile to analyse in more
detail how multi-step approaches e.g.\ by
    \citet{Glassermanetal05} have to be modified for successful application on the estimation
    of VaR
contributions as studied here.  
Further research on this issue could be useful.

%%%%%%%%%%%%%%
% References %
%%%%%%%%%%%%%%

%\bibliographystyle{plainnat}
%\bibliography{Literature}
%\input{Euler_allocation.bbl}

%\newpage
\setcounter{table}{0}
\refstepcounter{table}
%  \begin{table}[htp]
\begin{samepage}
    \begin{center}
\label{tab:1}
\parbox{15cm}{Table \thetable: \emph{Risk characteristics of assets. Identical for all
four sectors.}}\\[2ex]
\begin{tabular}{|c|c|c|c|c|c|}
\hline
Asset No. & PD & Exposure & LGD mean & LGD variance & Stand-alone VaR$_{99.9\%}$\\ \hline \hline
 1 & 0.02 & 25 & 0.5 & 0.125 & 24.846 \\ \hline
 2 & 0.02 & 25 & 0.5 & 0.125 & 24.846 \\ \hline
 3 & 0.02 & 25 & 0.5 & 0.03125 & 19.778 \\ \hline
 4 & 0.02 & 25 & 0.5 & 0.03125 & 19.778 \\ \hline
 5 & 0.02 & 5 & 0.5 & 0.125 & 4.969 \\ \hline
 6 & 0.02 & 5 & 0.5 & 0.125 & 4.969 \\ \hline
 7 & 0.02 & 5 & 0.5 & 0.03125 & 3.956 \\ \hline
 8 & 0.02 & 5 & 0.5 & 0.03125 & 3.956 \\ \hline
 9 & 0.02 & 1 & 0.5 & 0.125 & 0.994 \\ \hline
 10 & 0.02 & 1 & 0.5 & 0.125 & 0.994 \\ \hline
 11 & 0.02 & 1 & 0.5 & 0.03125 & 0.791 \\ \hline
 12 & 0.02 & 1 & 0.5 & 0.03125 & 0.791 \\ \hline
 13 & 0.005 & 25 & 0.5 & 0.125 & 22.613 \\ \hline
 14 & 0.005 & 25 & 0.5 & 0.125 & 22.613 \\ \hline
 15 & 0.005 & 25 & 0.5 & 0.03125 & 16.51 \\ \hline
 16 & 0.005 & 25 & 0.5 & 0.03125 & 16.51 \\ \hline
 17 & 0.005 & 5 & 0.5 & 0.125 & 4.523 \\ \hline
 18 & 0.005 & 5 & 0.5 & 0.125 & 4.523 \\ \hline
 19 & 0.005 & 5 & 0.5 & 0.03125 & 3.302 \\ \hline
 20 & 0.005 & 5 & 0.5 & 0.03125 & 3.302 \\ \hline
 21 & 0.005 & 1 & 0.5 & 0.125 & 0.905 \\ \hline
 22 & 0.005 & 1 & 0.5 & 0.125 & 0.905 \\ \hline
 23 & 0.005 & 1 & 0.5 & 0.03125 & 0.66 \\ \hline
 24 & 0.005 & 1 & 0.5 & 0.03125 & 0.66 \\ \hline
\end{tabular}
\end{center}
\end{samepage}

\refstepcounter{table}
%  \begin{table}[htp]
    \begin{center}
\label{tab:2}
\parbox{15cm}{Table \thetable: \emph{Standard Monte Carlo and importance sampling estimates
of portfolio VaR at 99.9\% level and coefficients of variation of estimates.}}\\[2ex]
\begin{tabular}{|l|c|c|c|}
\hline
Sampling method & 1st run & Mean of all runs & Coef.\ of var. \\ \hline \hline
Standard MC & 69.29 & 69.22 & 1.88\% \\ \hline
Imp.~samp.\ MC& 68.4 & 68.7 & 1.22\% \\ \hline
\end{tabular}
\end{center}

\refstepcounter{table}
%  \begin{table}[htp]
    \begin{center}
\label{tab:3}
\parbox{15cm}{Table \thetable: \emph{Standard Monte Carlo estimates
of conditional expectations of systematic factors conditional on ``Loss equals portfolio VaR at 99.9\% level''
and coefficients of variation of estimates.}}\\[2ex]
\begin{tabular}{|l|c|c|c|}
\hline
Factor & 1st run & Mean of all runs & Coef.\ of var. \\ \hline \hline
1 & -1.327 & -1.449 & 16.38\% \\ \hline
2 & -1.313 & -1.471 & 14.23\% \\ \hline
3 & -0.852 & -0.854 & 22.04\% \\ \hline
4 & -1.244 & -0.87 & 34.15\% \\ \hline
\end{tabular}
\end{center}

%\newpage

\refstepcounter{table}
%  \begin{table}[htp]
    \begin{center}
\label{tab:4}
\parbox{15cm}{Table \thetable: \emph{Portfolio loss distribution characteristics
mean of loss, standard deviation of loss, and probability of observing positive losses
for original and importance sampling distribution. Characteristics for original
distribution are calculated without simulation only once before the 1st simulation run.}}\\[2ex]
\begin{tabular}{|l|c||c|c|c|}
\hline
Sector & Original distribution & \multicolumn{3}{c|}{Importance sampling distribution}    \\\hline
 & & 1st run & Mean of all runs & Coef.\ of var. \\ \hline \hline
Expected loss & 6.2  & 70.101 & 73.204 & 15.73\% \\ \hline
Stddev of loss & 10.359  & 41.406 & 42.204 & 8.14\% \\ \hline
Prob.~pos.~loss & 59.411\% & 99.911\% & 99.881\% & 0.15\% \\ \hline
\end{tabular}
\end{center}

\refstepcounter{table}
%  \begin{table}[htp]
    \begin{center}
\label{tab:5}
\parbox{15cm}{Table \thetable: \emph{Standard Monte Carlo and importance sampling estimates
of stand-alone sector-VaR at 99.9\% level
and coefficients of variation of estimates.}}\\[2ex]
\begin{tabular}{|l|c|c|c||c|c|c|}
\hline
Sector & \multicolumn{3}{c||}{Standard MC} & \multicolumn{3}{c|}{Importance sampling}    \\\hline
 & 1st run & Mean of all runs & Coef.\ of var. & 1st run & Mean of all runs & Coef.\ of var. \\ \hline \hline
1 & 39.434 & 41.379 & 2.79\% & 40.366 & 41.253 & 1.95\% \\ \hline
2 & 42.302 & 41.524 & 2.67\% & 39.969 & 41.487 & 2.53\% \\ \hline
3 & 40.954 & 41.504 & 2.49\% & 40.369 & 41.211 & 3.37\% \\ \hline
4 & 41.15 & 41.285 & 2.08\% & 42.674 & 41.623 & 4.71\% \\ \hline
\end{tabular}
\end{center}

\refstepcounter{table}
%  \begin{table}[htp]
    \begin{center}
\label{tab:6}
\parbox{15cm}{Table \thetable: \emph{Ratio of sum of asset VaR contributions and portfolio-wide VaR (at 99.9\% level) for Standard Monte Carlo, importance sampling, importance sampling with 50\% of bandwidth, and
importance sampling with 200\% of bandwidth.}}\\[2ex]
\begin{tabular}{|l|c|c|}
\hline
Sampling method & 1st run & With mean contributions and VaRs of all runs \\ \hline \hline
Standard MC & 99.69\% & 99.6\% \\ \hline
Imp.~samp.~MC   & 96.84\% & 96.2\% \\ \hline
50\% bandwidth & 99.15\% & 98.96\% \\ \hline
200\% bandwidth & 85.31\% & 85.29\% \\ \hline
\end{tabular}
\end{center}

\newpage
\begin{samepage}
\refstepcounter{table}
%  \begin{table}[htp]
    \begin{center}
\small
\label{tab:7}
\parbox{15cm}{Table \thetable: \emph{Standard Monte Carlo and importance sampling estimates
of asset contributions to portfolio-wide VaR at 99.9\% level
and coefficients of variation of estimates. Only sectors 1 and 3. See Figure \ref{fig:1} for a 
graphical comparison of the coefficients of variation.}}\\[2ex]
\begin{tabular}{|l|l|c|c|c||c|c|c|}
\hline
Sector & Asset & \multicolumn{3}{c||}{Standard MC} & \multicolumn{3}{c|}{Importance sampling}    \\\hline
 & & 1st run & Mean of all runs & Coef.\ of var. & 1st run & Mean of all runs & Coef.\ of var. \\ \hline \hline
1 & 1 & 3.73991 & 4.13472 & 39.57\% & 3.5955 & 4.1292 & 28.43\% \\ \hline
1 & 2 & 2.48008 & 3.98558 & 42.72\% & 4.57874 & 4.10322 & 14.51\% \\ \hline
1 & 3 & 1.86887 & 2.78517 & 43.7\% & 2.13936 & 2.70369 & 16.58\% \\ \hline
1 & 4 & 1.67832 & 2.98133 & 41.23\% & 3.90476 & 2.79776 & 23.62\% \\ \hline
1 & 5 & 0.00852 & 0.32861 & 64.19\% & 0.29531 & 0.261 & 25.11\% \\ \hline
1 & 6 & 0.25436 & 0.26363 & 68.58\% & 0.36958 & 0.28069 & 20.92\% \\ \hline
1 & 7 & 0.01275 & 0.29346 & 60.74\% & 0.28293 & 0.23594 & 19.17\% \\ \hline
1 & 8 & 0.10971 & 0.24488 & 61.28\% & 0.2703 & 0.24198 & 20.01\% \\ \hline
1 & 9 & 0.04055 & 0.04092 & 64.45\% & 0.05336 & 0.04391 & 26.94\% \\ \hline
1 & 10 & 0.03549 & 0.03701 & 75.96\% & 0.04931 & 0.03976 & 14.91\% \\ \hline
1 & 11 & 0.03 & 0.04138 & 62.7\% & 0.0448 & 0.04162 & 16.32\% \\ \hline
1 & 12 & 0.05875 & 0.05802 & 66.77\% & 0.0462 & 0.04339 & 22.3\% \\ \hline
1 & 13 & 0.90201 & 1.26936 & 98.54\% & 1.35922 & 1.33622 & 30.04\% \\ \hline
1 & 14 & 0.72223 & 0.99663 & 73.19\% & 1.96059 & 1.27473 & 22.41\% \\ \hline
1 & 15 & 0.00017 & 0.58035 & 85.56\% & 0.81998 & 0.97163 & 28.52\% \\ \hline
1 & 16 & 0.88218 & 0.78081 & 65.45\% & 1.191 & 0.87679 & 20.11\% \\ \hline
1 & 17 & 0.01785 & 0.09645 & 126.93\% & 0.08795 & 0.09451 & 58.51\% \\ \hline
1 & 18 & 0.02523 & 0.10587 & 100.98\% & 0.05491 & 0.08176 & 24.24\% \\ \hline
1 & 19 & 0.00297 & 0.08758 & 109.65\% & 0.09981 & 0.08263 & 34.23\% \\ \hline
1 & 20 & 0.18965 & 0.10518 & 120.03\% & 0.07885 & 0.08013 & 39.55\% \\ \hline
1 & 21 & 0.01695 & 0.02205 & 94.88\% & 0.02105 & 0.01382 & 31.04\% \\ \hline
1 & 22 & 0.00065 & 0.0101 & 182.79\% & 0.01315 & 0.01494 & 41.93\% \\ \hline
1 & 23 & 0.01227 & 0.01647 & 115.13\% & 0.01212 & 0.01434 & 22.03\% \\ \hline
1 & 24 & 0.00011 & 0.00962 & 123.95\% & 0.00676 & 0.01353 & 25.55\% \\ \hline
3 & 1 & 4.22536 & 2.87626 & 50.72\% & 3.20845 & 2.93319 & 23.72\% \\ \hline
3 & 2 & 4.16669 & 3.57096 & 35.82\% & 2.42476 & 2.82777 & 22.43\% \\ \hline
3 & 3 & 1.85084 & 1.82327 & 60.66\% & 1.62841 & 1.94189 & 26.81\% \\ \hline
3 & 4 & 1.86416 & 1.71489 & 44.15\% & 1.67836 & 1.92037 & 28.6\% \\ \hline
3 & 5 & 0.00607 & 0.15006 & 86.14\% & 0.30478 & 0.17055 & 31.19\% \\ \hline
3 & 6 & 0.26857 & 0.23839 & 69.98\% & 0.29962 & 0.15907 & 29.65\% \\ \hline
3 & 7 & 0.18092 & 0.16485 & 75.66\% & 0.13607 & 0.17009 & 27.24\% \\ \hline
3 & 8 & 0.31751 & 0.1394 & 76.11\% & 0.23939 & 0.15737 & 35.18\% \\ \hline
3 & 9 & 0.03608 & 0.02836 & 91.3\% & 0.0307 & 0.02636 & 30.25\% \\ \hline
3 & 10 & 0.04145 & 0.03967 & 76.9\% & 0.02157 & 0.02805 & 33.13\% \\ \hline
3 & 11 & 0.01261 & 0.02728 & 93.85\% & 0.03185 & 0.02711 & 35.09\% \\ \hline
3 & 12 & 0.01994 & 0.02723 & 85.55\% & 0.01969 & 0.02961 & 23.47\% \\ \hline
3 & 13 & 1.38858 & 1.22577 & 88.2\% & 0.91357 & 1.04637 & 85.95\% \\ \hline
3 & 14 & 1.68561 & 1.0336 & 77.07\% & 0.80861 & 0.87244 & 29.14\% \\ \hline
3 & 15 & 0.64665 & 0.54598 & 93.05\% & 0.53973 & 0.58639 & 29.58\% \\ \hline
3 & 16 & 1.06298 & 0.80942 & 75.61\% & 0.73897 & 0.62549 & 42.79\% \\ \hline
3 & 17 & 0.00001 & 0.04572 & 201.77\% & 0.07141 & 0.05265 & 46.15\% \\ \hline
3 & 18 & 0.00037 & 0.04821 & 178.53\% & 0.06335 & 0.0534 & 61.58\% \\ \hline
3 & 19 & 0.00561 & 0.06153 & 126.15\% & 0.05904 & 0.04706 & 36.99\% \\ \hline
3 & 20 & 0.00008 & 0.03582 & 165.25\% & 0.04812 & 0.05005 & 80.54\% \\ \hline
3 & 21 & 0.02334 & 0.016 & 132.45\% & 0.00505 & 0.0107 & 68.05\% \\ \hline
3 & 22 & 0.04998 & 0.00533 & 230.52\% & 0.01307 & 0.00966 & 60.67\% \\ \hline
3 & 23 & 0 & 0.01009 & 192.67\% & 0.00971 & 0.00835 & 65.11\% \\ \hline
3 & 24 & 0 & 0.00651 & 152.59\% & 0.0074 & 0.01027 & 58.74\% \\ \hline
\end{tabular}
\end{center}
\end{samepage}

\newpage
\setcounter{figure}{0}
\refstepcounter{figure}
\begin{figure}[htb]
  \begin{center}
%\centering
  \parbox{14.0cm}{Figure \thefigure:
  \emph{Graphical comparison of coefficients of variation from Table \ref{tab:7}.}}
\label{fig:1}\\[1ex]
\begin{turn}{270}
\resizebox{\height}{14.0cm}{\includegraphics[width=14.0cm]{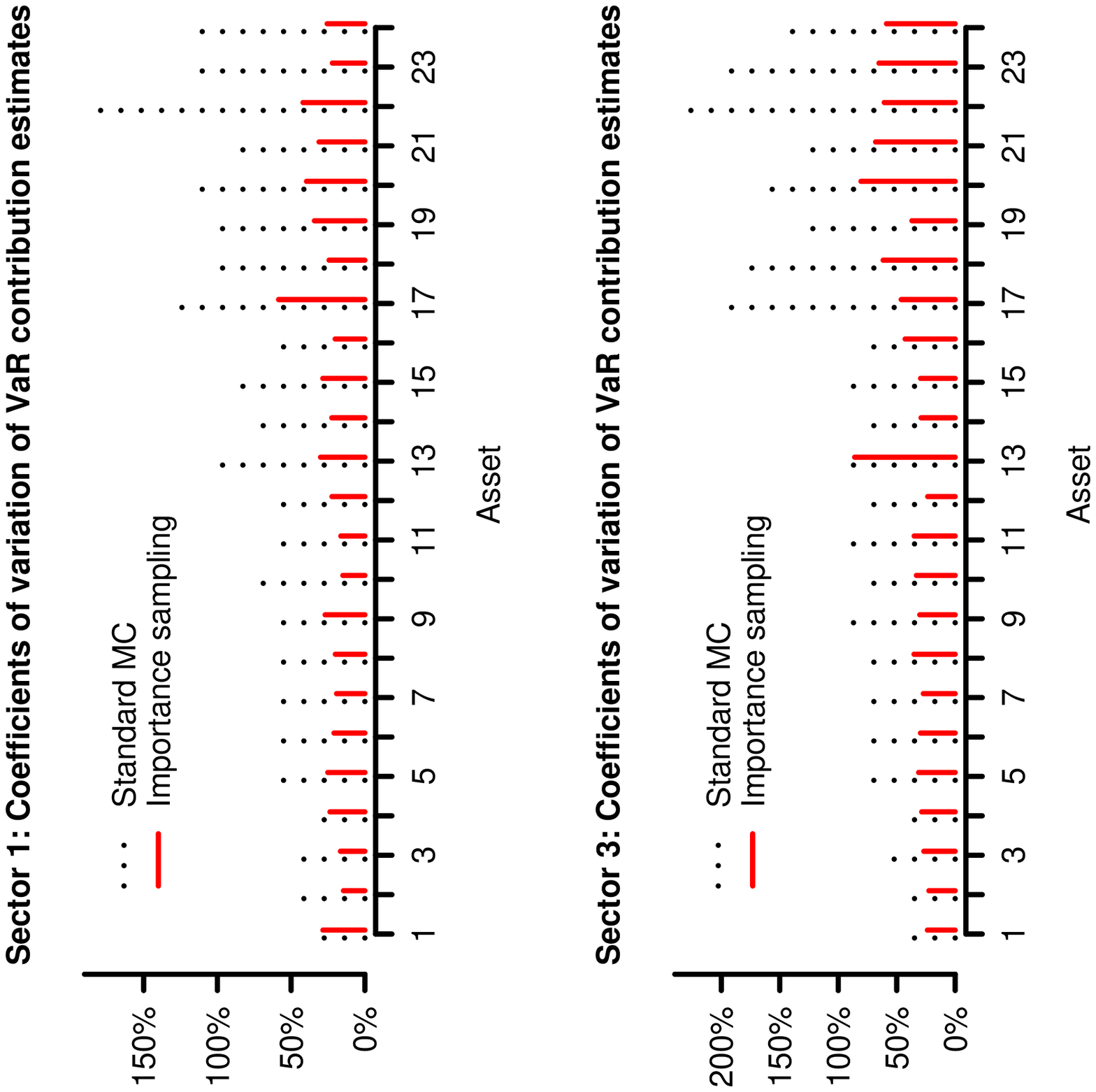}}
\end{turn}
\end{center}
\end{figure}

\newpage
\begin{samepage}
\refstepcounter{table}
%  \begin{table}[htp]
    \begin{center}
\small
\label{tab:8}
\parbox{15cm}{Table \thetable: \emph{Importance sampling estimates with half and double
of bandwidth according to \eqref{eq:Silverman}
for asset contributions to portfolio-wide VaR at 99.9\% level
and coefficients of variation of estimates. Only sectors 1 and 3.}}\\[2ex]
\begin{tabular}{|l|l|c|c|c||c|c|c|}
\hline
Sector & Asset & \multicolumn{3}{c||}{Imp.~samp., 50\% bandwidth} & \multicolumn{3}{c|}{Imp.~samp., 200\% bandwidth}    \\\hline
 & & 1st run & Mean of all runs & Coef.\ of var. & 1st run & Mean of all runs & Coef.\ of var. \\ \hline \hline
1 & 1 & 3.26946 & 4.29182 & 38.02\% & 4.05739 & 4.00807 & 16.07\% \\ \hline
1 & 2 & 4.66724 & 3.94153 & 16.02\% & 4.03265 & 4.09244 & 12.27\% \\ \hline
1 & 3 & 1.98064 & 2.64532 & 19.71\% & 2.14906 & 2.72078 & 16.23\% \\ \hline
1 & 4 & 3.80636 & 2.57832 & 23.22\% & 3.34566 & 2.8715 & 18.5\% \\ \hline
1 & 5 & 0.26897 & 0.26431 & 41.28\% & 0.28143 & 0.25674 & 16.95\% \\ \hline
1 & 6 & 0.33207 & 0.29524 & 34.96\% & 0.30172 & 0.26413 & 12.17\% \\ \hline
1 & 7 & 0.26506 & 0.23521 & 19.71\% & 0.24753 & 0.2373 & 16.24\% \\ \hline
1 & 8 & 0.24949 & 0.24084 & 25.59\% & 0.23923 & 0.24002 & 14.17\% \\ \hline
1 & 9 & 0.05815 & 0.04473 & 34.66\% & 0.04556 & 0.04185 & 15.54\% \\ \hline
1 & 10 & 0.05494 & 0.03912 & 21.92\% & 0.04404 & 0.04157 & 10.76\% \\ \hline
1 & 11 & 0.04378 & 0.04123 & 20.37\% & 0.04292 & 0.04141 & 12.44\% \\ \hline
1 & 12 & 0.04586 & 0.04448 & 38.31\% & 0.04482 & 0.04233 & 12.57\% \\ \hline
1 & 13 & 1.34413 & 1.3799 & 43.47\% & 1.17617 & 1.26651 & 18.33\% \\ \hline
1 & 14 & 1.90522 & 1.31631 & 31.05\% & 1.56665 & 1.21016 & 17.01\% \\ \hline
1 & 15 & 0.84645 & 1.04631 & 45.4\% & 0.85945 & 0.88251 & 20.24\% \\ \hline
1 & 16 & 1.27248 & 0.87883 & 36.64\% & 0.92764 & 0.88677 & 18.53\% \\ \hline
1 & 17 & 0.05814 & 0.09153 & 57.75\% & 0.10907 & 0.08864 & 43.24\% \\ \hline
1 & 18 & 0.04969 & 0.08279 & 30.12\% & 0.06831 & 0.0812 & 17.81\% \\ \hline
1 & 19 & 0.10725 & 0.08509 & 52.91\% & 0.08653 & 0.07897 & 19.23\% \\ \hline
1 & 20 & 0.0893 & 0.08236 & 62.2\% & 0.06944 & 0.08108 & 21.85\% \\ \hline
1 & 21 & 0.02114 & 0.01322 & 47.88\% & 0.01951 & 0.01459 & 30.43\% \\ \hline
1 & 22 & 0.01281 & 0.01556 & 50.13\% & 0.01248 & 0.01424 & 25.29\% \\ \hline
1 & 23 & 0.01032 & 0.01462 & 31.32\% & 0.01293 & 0.01341 & 18.77\% \\ \hline
1 & 24 & 0.0058 & 0.01343 & 38.54\% & 0.00996 & 0.0133 & 20.48\% \\ \hline
3 & 1 & 4.01468 & 3.04237 & 34.49\% & 2.93766 & 3.06196 & 16.15\% \\ \hline
3 & 2 & 2.43925 & 2.88584 & 33.77\% & 2.67336 & 2.9815 & 12.16\% \\ \hline
3 & 3 & 1.66845 & 1.96384 & 36.29\% & 1.94701 & 2.05924 & 24.65\% \\ \hline
3 & 4 & 1.61369 & 1.89573 & 45.73\% & 1.79944 & 2.00395 & 17.04\% \\ \hline
3 & 5 & 0.39934 & 0.17245 & 46.26\% & 0.22462 & 0.17462 & 19.1\% \\ \hline
3 & 6 & 0.30854 & 0.15244 & 34.22\% & 0.26823 & 0.17591 & 19.79\% \\ \hline
3 & 7 & 0.14593 & 0.17199 & 35.86\% & 0.17599 & 0.1761 & 19.23\% \\ \hline
3 & 8 & 0.3268 & 0.15174 & 46.17\% & 0.22091 & 0.17559 & 29.68\% \\ \hline
3 & 9 & 0.03699 & 0.02606 & 45.78\% & 0.03159 & 0.02761 & 21.52\% \\ \hline
3 & 10 & 0.0174 & 0.02805 & 46.25\% & 0.0248 & 0.02895 & 22.52\% \\ \hline
3 & 11 & 0.03086 & 0.02792 & 47.53\% & 0.0313 & 0.02822 & 21.98\% \\ \hline
3 & 12 & 0.01679 & 0.02996 & 35.62\% & 0.02146 & 0.02889 & 15.67\% \\ \hline
3 & 13 & 0.92607 & 1.08528 & 113.38\% & 0.88461 & 1.01339 & 44.1\% \\ \hline
3 & 14 & 0.66095 & 0.86246 & 41.28\% & 0.89759 & 0.93943 & 25.98\% \\ \hline
3 & 15 & 0.32995 & 0.57864 & 39.28\% & 0.71209 & 0.63163 & 17.53\% \\ \hline
3 & 16 & 0.76557 & 0.59353 & 66.16\% & 0.80199 & 0.63174 & 27.11\% \\ \hline
3 & 17 & 0.10155 & 0.05401 & 50.01\% & 0.05136 & 0.05413 & 33.84\% \\ \hline
3 & 18 & 0.08667 & 0.05428 & 85.39\% & 0.06688 & 0.05492 & 35.01\% \\ \hline
3 & 19 & 0.05978 & 0.0426 & 57.18\% & 0.05109 & 0.04891 & 26.68\% \\ \hline
3 & 20 & 0.03033 & 0.05615 & 121.29\% & 0.06657 & 0.0519 & 34.31\% \\ \hline
3 & 21 & 0.00573 & 0.01152 & 102.81\% & 0.00803 & 0.01021 & 44.1\% \\ \hline
3 & 22 & 0.00939 & 0.00933 & 100.7\% & 0.01451 & 0.01001 & 43.49\% \\ \hline
3 & 23 & 0.01156 & 0.00858 & 102.21\% & 0.00945 & 0.00872 & 37\% \\ \hline
3 & 24 & 0.00558 & 0.0119 & 83.98\% & 0.00905 & 0.00922 & 34.99\% \\ \hline
\end{tabular}
\end{center}
\end{samepage}

\newpage

\refstepcounter{table}
%  \begin{table}[htp]
    \begin{center}
\label{tab:9}
\parbox{15cm}{Table \thetable: \emph{Standard Monte Carlo and importance sampling estimates
of sector contributions to VaR at 99.9\% level
and coefficients of variation of estimates.}}\\[2ex]
\begin{tabular}{|l|c|c|c||c|c|c|}
\hline
Sector & \multicolumn{3}{c||}{Standard MC} & \multicolumn{3}{c|}{Importance sampling}    \\\hline
 & 1st run & Mean of all runs & Coef.\ of var. & 1st run & Mean of all runs & Coef.\ of var. \\ \hline \hline
1 & 13.09 & 19.275 & 23.61\% & 21.336 & 19.777 & 7.92\% \\ \hline
2 & 20.507 & 20.512 & 16.39\% & 18.51 & 19.683 & 10.29\% \\ \hline
3 & 17.853 & 14.645 & 23.53\% & 13.302 & 13.764 & 14.01\% \\ \hline
4 & 17.845 & 14.784 & 28.72\% & 15.257 & 15.474 & 15.31\% \\ \hline
\end{tabular}
\end{center}

\refstepcounter{table}
%  \begin{table}[htp]
    \begin{center}
\label{tab:10}
\parbox{15cm}{Table \thetable: \emph{Standard Monte Carlo and importance sampling estimates
of portfolio-wide diversification index with respect to VaR at 99.9\% level, in asset-level and
sector-level context,
and coefficients of variation of estimates.}}\\[2ex]
\begin{tabular}{|l|c|c|c||c|c|c|}
\hline
Context & \multicolumn{3}{c||}{Standard MC} & \multicolumn{3}{c|}{Importance sampling}    \\\hline
 & 1st run & Mean of all runs & Coef.\ of var. & 1st run & Mean of all runs & Coef.\ of var. \\ \hline \hline
Assets & 7.65\% & 7.64\% & 2.06\% & 7.54\% & 7.58\% & 1.34\% \\ \hline
Sectors & 40.02\% & 39.51\% & 2.09\% & 39.58\% & 39.21\% & 2.03\% \\ \hline
\end{tabular}
\end{center}

\newpage
\begin{samepage}
\refstepcounter{table}
%  \begin{table}[htp]
    \begin{center}
\small
\label{tab:11}
\parbox{15cm}{Table \thetable: \emph{Standard Monte Carlo and importance sampling estimates
of marginal diversification indices at asset-level with respect to portfolio-wide VaR at 99.9\% level
and coefficients of variation of estimates. Only sectors 1 and 3.}}\\[2ex]
\begin{tabular}{|l|l|c|c|c||c|c|c|}
\hline
Sector & Asset & \multicolumn{3}{c||}{Standard MC} & \multicolumn{3}{c|}{Importance sampling} \\\hline
 & & 1st run & Mean of all runs & Coef.\ of var. & 1st run & Mean of all runs & Coef.\ of var. \\ \hline \hline
1 & 1 & 14.19\% & 15.8\% & 42.12\% & 13.61\% & 15.79\% & 30.35\% \\ \hline
1 & 2 & 9.07\% & 15.19\% & 45.6\% & 17.62\% & 15.69\% & 15.5\% \\ \hline
1 & 3 & 8.29\% & 12.98\% & 48.03\% & 9.66\% & 12.56\% & 18.33\% \\ \hline
1 & 4 & 7.31\% & 13.99\% & 45.01\% & 18.74\% & 13.05\% & 26.05\% \\ \hline
1 & 5 & -0.85\% & 5.66\% & 75.78\% & 4.97\% & 4.27\% & 31.38\% \\ \hline
1 & 6 & 4.15\% & 4.34\% & 84.72\% & 6.48\% & 4.67\% & 25.68\% \\ \hline
1 & 7 & -0.96\% & 6.23\% & 73.27\% & 5.94\% & 4.73\% & 24.63\% \\ \hline
1 & 8 & 1.53\% & 4.99\% & 77.09\% & 5.62\% & 4.88\% & 25.47\% \\ \hline
1 & 9 & 3.1\% & 3.14\% & 85.41\% & 4.39\% & 3.42\% & 35.3\% \\ \hline
1 & 10 & 2.59\% & 2.74\% & 104.24\% & 3.98\% & 3\% & 20.21\% \\ \hline
1 & 11 & 2.56\% & 4.01\% & 82.79\% & 4.43\% & 4.01\% & 21.77\% \\ \hline
1 & 12 & 6.24\% & 6.15\% & 80.74\% & 4.61\% & 4.24\% & 29.35\% \\ \hline
1 & 13 & 3.72\% & 5.35\% & 103.64\% & 5.76\% & 5.66\% & 31.61\% \\ \hline
1 & 14 & 2.93\% & 4.14\% & 78.11\% & 8.44\% & 5.39\% & 23.6\% \\ \hline
1 & 15 & -0.38\% & 3.15\% & 95.93\% & 4.61\% & 5.53\% & 30.56\% \\ \hline
1 & 16 & 4.98\% & 4.37\% & 71.17\% & 6.87\% & 4.96\% & 21.72\% \\ \hline
1 & 17 & 0.12\% & 1.86\% & 145.91\% & 1.67\% & 1.81\% & 67.75\% \\ \hline
1 & 18 & 0.28\% & 2.07\% & 114.54\% & 0.93\% & 1.53\% & 28.81\% \\ \hline
1 & 19 & -0.29\% & 2.28\% & 127.99\% & 2.65\% & 2.13\% & 40.63\% \\ \hline
1 & 20 & 5.39\% & 2.82\% & 136.3\% & 2.01\% & 2.05\% & 47.16\% \\ \hline
1 & 21 & 1.6\% & 2.17\% & 107.04\% & 2.05\% & 1.25\% & 38.24\% \\ \hline
1 & 22 & -0.21\% & 0.84\% & 243.25\% & 1.18\% & 1.37\% & 50.75\% \\ \hline
1 & 23 & 1.48\% & 2.12\% & 135.81\% & 1.45\% & 1.79\% & 26.88\% \\ \hline
1 & 24 & -0.36\% & 1.08\% & 167.71\% & 0.64\% & 1.67\% & 31.62\% \\ \hline
3 & 1 & 16.16\% & 10.68\% & 55.57\% & 12.03\% & 10.91\% & 26.02\% \\ \hline
3 & 2 & 15.93\% & 13.5\% & 38.53\% & 8.84\% & 10.48\% & 24.7\% \\ \hline
3 & 3 & 8.2\% & 8.05\% & 70.34\% & 7.04\% & 8.65\% & 30.96\% \\ \hline
3 & 4 & 8.26\% & 7.5\% & 51.71\% & 7.3\% & 8.54\% & 33.08\% \\ \hline
3 & 5 & -0.9\% & 2.03\% & 129.42\% & 5.16\% & 2.42\% & 44.94\% \\ \hline
3 & 6 & 4.44\% & 3.83\% & 88.67\% & 5.06\% & 2.19\% & 44.09\% \\ \hline
3 & 7 & 3.35\% & 2.94\% & 108.77\% & 2.17\% & 3.04\% & 39.1\% \\ \hline
3 & 8 & 6.85\% & 2.28\% & 118.93\% & 4.82\% & 2.71\% & 52.58\% \\ \hline
3 & 9 & 2.65\% & 1.86\% & 141.32\% & 2.08\% & 1.63\% & 49.99\% \\ \hline
3 & 10 & 3.19\% & 3.01\% & 102.98\% & 1.15\% & 1.8\% & 52.75\% \\ \hline
3 & 11 & 0.33\% & 2.21\% & 148.43\% & 2.77\% & 2.15\% & 56.96\% \\ \hline
3 & 12 & 1.27\% & 2.2\% & 135.55\% & 1.2\% & 2.47\% & 36.1\% \\ \hline
3 & 13 & 5.88\% & 5.16\% & 92.98\% & 3.78\% & 4.37\% & 91.65\% \\ \hline
3 & 14 & 7.2\% & 4.31\% & 82.07\% & 3.31\% & 3.59\% & 31.49\% \\ \hline
3 & 15 & 3.55\% & 2.94\% & 105.14\% & 2.9\% & 3.18\% & 33.27\% \\ \hline
3 & 16 & 6.08\% & 4.54\% & 81.95\% & 4.11\% & 3.42\% & 47.75\% \\ \hline
3 & 17 & -0.28\% & 0.74\% & 278.06\% & 1.3\% & 0.88\% & 61.29\% \\ \hline
3 & 18 & -0.27\% & 0.79\% & 241.35\% & 1.12\% & 0.9\% & 81.37\% \\ \hline
3 & 19 & -0.21\% & 1.49\% & 158.46\% & 1.41\% & 1.04\% & 51.06\% \\ \hline
3 & 20 & -0.38\% & 0.71\% & 254.37\% & 1.07\% & 1.13\% & 108.81\% \\ \hline
3 & 21 & 2.31\% & 1.5\% & 157.1\% & 0.27\% & 0.9\% & 89.92\% \\ \hline
3 & 22 & 5.26\% & 0.31\% & 435.57\% & 1.17\% & 0.79\% & 83.06\% \\ \hline
3 & 23 & -0.38\% & 1.15\% & 256.44\% & 1.09\% & 0.88\% & 94.5\% \\ \hline
3 & 24 & -0.38\% & 0.61\% & 248.39\% & 0.74\% & 1.17\% & 78.71\% \\ \hline
\end{tabular}
\end{center}
\end{samepage}

\newpage

\refstepcounter{table}
%  \begin{table}[htp]
    \begin{center}
\label{tab:12}
\parbox{15cm}{Table \thetable: \emph{Standard Monte Carlo and importance sampling estimates
of diversification indices at sector-level with respect to portfolio-wide VaR at 99.9\% level
and coefficients of variation of estimates.}}\\[2ex]
\begin{tabular}{|l|c|c|c||c|c|c|}
\hline
Sector & \multicolumn{3}{c||}{Standard MC} & \multicolumn{3}{c|}{Importance sampling}    \\\hline
 & 1st run & Mean of all runs & Coef.\ of var. & 1st run & Mean of all runs & Coef.\ of var. \\ \hline \hline
1 & 30.46\% & 44.51\% & 25.27\% & 51.01\% & 45.94\% & 9.39\% \\ \hline
2 & 46.52\% & 47.44\% & 17.04\% & 44.16\% & 45.43\% & 10.78\% \\ \hline
3 & 41.37\% & 32.77\% & 25.98\% & 30.24\% & 30.76\% & 15.13\% \\ \hline
4 & 41.15\% & 33.3\% & 32.32\% & 33.32\% & 34.73\% & 15.97\% \\ \hline
\end{tabular}
\end{center}

\end{document}